\documentclass[12pt]{article}
\usepackage{amsmath}
\usepackage{amsfonts}
\usepackage{amsthm}
\usepackage{graphicx}
\usepackage{hyperref}

\numberwithin{equation}{section} 

\newtheorem{thm}{Theorem}[section]

\newtheorem{prop}[thm]{Proposition}

\newtheorem{defn}[thm]{Definition}

\newcommand{\N}{\mathbb{N}}         
\newcommand{\Wnn}{W_{n,n}}  
\newcommand{\Knn}{K_{n,n}}  
\newcommand{\tfp}{\tilde{F}_\Phi} 

\begin{document}

\title{Experimental analysis of lattice walks}

\author{
Anthony Zaleski\thanks{Department of Mathematics, Rutgers University (New Brunswick), 110 Frelinghuysen Road, Piscataway, NJ 08854-8019, USA.}
}

\maketitle
\begin{abstract}
The book \cite{feller} by Feller discusses statistics corresponding to sequences of coin tosses, with a dollar being won or lost depending on the outcome of each toss.  This is equivalent to analyzing walks in the plane with each step being one unit up or right.  

One statistic corresponding to a walk is the ``number of losing times," i.e., the number of times the player has negative money, or the number of times the 2-D walk is below the line $y=x$.  Two main results for this random variable are that it follows the ``discrete arcsine law" when the sample space is walks of length $n$, and it is described by the Chung-Feller generating function when the sample space is walks to $(n,n)$.

In \cite{zeilberger} and the accompanying Maple package, Zeilberger outlines a way to obtain Feller's results with the help of a computer.  In fact, Zeilberger's Maple package computes the ``grand generating function," which, in a single blow, captures information about \emph{all} of the walk statistics discussed in \cite{feller}.

In this paper, we continue to investigate walks using computer methods.  However, we shall introduce an approach different from that of Zeilberger in \cite{zeilberger}.  In \cite{zeilberger}, computer algebra is used to exactly compute the weight enumerator over \emph{all} walks---an \emph{infinite} sum expressed as an algebraic function.  Our procedures input a \emph{numeric} value of $n$ and use recursion to find the weight enumerator over (for example) all walks of length $n$---a \emph{finite} polynomial.  Then, by generating these polynomials for many values of $n$, we are able to guess behavior of the moments of certain statistics as the walk length tends to $\infty$.

The advantage of this method is that it is easily applied to more general problems not amenable to an analytic approach.  For example, we use it to analyze walks in three dimensions and walks where more general steps (e.g., diagonal steps) are allowed.  
\end{abstract}
\section{Introduction}\label{sec:intro}
\subsection{Coin tossing and walks in the plane}
Suppose a gambler tosses a coin finitely many times, winning a dollar whenever heads comes up and losing a dollar when tails appears.  The evolution of the game can be described by a string $w=w_1w_2\cdots w_n$, where $w_i \in \{-1,1\}$ describes the outcome of the $i^{th}$ coin toss. 

By making the association $-1=r$ (a right step) and $1=u$ (an up step), we can also interpret a string $w$ as a walk in $\N^2$ starting from $(0,0)$.  For example, the string $ur$ represents a game in which a dollar is won and then lost.  Equivalently, $ur$ is a walk from $(0,0)\to(1,0)\to(1,1)$.  We shall use the terms ``game," ``string," and ``walk" interchangeably. 

\begin{defn}\label{def:stat}
Let $W$ be the set of all such walks.  We have the following walk statistics (functions from $W \to \N$):
\begin{itemize}
\item the length (number of steps), $l(w)$;
\item the number of losing times (points where the walk is below $y=x$): 
$$
a_1(w):=\left| \left\{i : \sum_{j=1}^i w_j <0 \text{ or } \left( \sum_{j=1}^i w_j =0 \text{ and } \sum_{j=1}^{i-1} w_j <0\right) \right\}\right|;
$$
\item the number of break-even times (points on $y=x$):
$$
a_2(w):=\left| \left\{ i : \sum_{j=1}^i w_i = 0\right\} \right|;
$$
\item the last break-even time:
$$
a_3(w):=\max \left\{ i : \sum_{j=1}^i w_j = 0 \text{ and }\sum_{j=1}^i w_j > 0 \text{ for } r>i \right\};
$$
\item and the number of sign-changes (points where the walk crosses $y=x$):
$$
a_4(w):=\left| \left\{ i : \sum_{j=1}^{i-1} w_j \cdot \sum_{j=1}^{i+1} w_j < 0 \right\} \right|.
$$
\end{itemize}

Further, we define $W_n$ to be the set of $n$-step walks in $W$, and $W_{n,n}$ to be the walks to $(n,n)$.
\end{defn}

\subsection{Old results}
Theorem III.4.1 in \cite{feller} states that there are ${2k \choose k} {2n-2k \choose n-k}$ walks $w\in W_n$ satisfying $a_1(w)=2k$. For $n$ large and fixed, the distribution of $a_1$ resembles (modulo scaling) $1/\sqrt{x(1-x)}$, so it is sometimes called the \emph{discrete arcsine} distribution. It is $u$-shaped, meaning that, surprisingly (or not surprisingly, if you believe in luck), most walks are either winning for most of their duration or losing on the majority of flips.

Theorem III.9 in \cite{feller}, the Chung-Feller theorem, says that the number of walks to $(n,n)$ with $2k$ losing times is given by a Catalan number, and is independent of $k$. In terms of generating functions:
$$
\sum_{w \in \Wnn} t^{a_1}(w) =\frac{1}{n-1} {2n \choose n} \sum_{k=0}^n t^{2k}.
$$

Finally, in \cite{zeilberger}, Zeilberger uses Maple to evaluate the ``grand generating function"
$$
\sum_{w \in W} z^{l(w)} t_1^{a_1(w)}t_2^{a_2(w)}t_3^{a_3(w)}t_4^{a_4(w)}
$$
as a (very messy) algebraic function of $z,t_1,t_2,t_3,t_4$.

These results are illuminating for this problem, but they are gotten through ad hoc methods.  So, for example, it is not obvious how to derive analogous results for walks in higher dimensions, or walks where non-standard steps are allowed. 

Here, we shall give an alternative method to analyze the statistics of a very general class of walks and approximate the long-run behavior of their moments.  But first, let us see if we can discover a few more exact results using the analytical approach.
\section{Exact results}
\subsection{Moments of up-right walks to $(n,n)$}
Suppose we uniformly randomly pick a walk $w\in \Wnn$.  Then we can think of $a_1(w)$ as a random variable.  For each $n$, $a_1 | \Wnn$ ($|$ = ``with sample space") has a certain distribution, so it is natural to wonder about the limiting distribution as $n\to \infty$.  For example, is it asymptotically normal?

Recall that the moments of a random variable are a fingerprint of its distribution; for example, the (standardized central) moments $0,1,0,3,0,15,0,105,\dots$ uniquely determine the normal distribution.  

In this case of $a_1 | \Wnn$, we can find the moments in terms of $n$, which is not surprising since $a_1$ is essentially uniform by the Chung-Feller rule.

Using the procedure \verb+ChungFeller+ in \cite{zeilberger}, we find
$$
\sum_{w \in \Wnn, \, n \in \N} z^n t^{a_1(w)}=\frac{2}{ \sqrt {-4\,z+1}+\sqrt {-4\,z{t}^{2}+1}}.
$$

Now we use \verb+convert(%,FPS,z)+ to convert this function to a formal power series in $z$.  By looking at the coefficient of $z^n$, we obtain the generating function
$$
F_n(t)=\sum_{w \in \Wnn} t^{a_1(w)},
$$
\emph{as a function of $n$!} To see it for yourself, use \verb+ChungFellerGF(t,n)+.

Now it is easy to compute the moments as functions of $n$, which is done in \verb+ChungFellerMoment+.  You can easily verify the following:
\begin{prop}\label{prop:cfm}
The number of losing times of a walk chosen uniformly randomly from $\Wnn$ has mean $n$ and variance $n^2/3+2n/3$, and its third  through tenth standardized moments about the mean approach $0,9/5,0,27/7,0,$ $9,0,243/11$ as $n \to \infty$.
\end{prop}

Analogously, for $a_2$, the number of visits to the diagonal $y=x$, we have
$$
F(z,t):=\sum_{w \in \Wnn, \, n \in \N} z^n t^{a_2(w)}=\frac{1}{t\sqrt {-4\,z+1}-t+1}.
$$

Unfortunately, Maple cannot convert this to a formal power series in $z$.  However, $F_t(z,1)$ \emph{is} convertible to a formal power series, so we can compute $\mathbb{E}[a_2 | \Wnn]=[z^n]F_t(z,1)/{2n \choose n}$, as a function of $n$.  In a similar way, we can find higher moments: the idea is to repeatedly apply the operator $t\partial_t$, substitute $t=1$, and \emph{then} expand as a formal power series in $z$.  (Some conversions are necessary to get the central moments.)  

The moments of $a_2$ are surprisingly complicated in comparison with those of $a_1$:
\begin{prop}
The number of visits to $y=x$ of a walk chosen uniformly randomly from $\Wnn$ has mean and variance
$$
{\frac {- \left( 2\,n \right) !+{4}^{n} \left( n! \right) ^{2}}{
 \left( 2\,n \right) !}},
 -{\frac {{16}^{n} \left( n! \right) ^{4}+{4}^{n} \left( n! \right) ^{2
} \left( 2\,n \right) !-4\,n \left(  \left( 2\,n \right) ! \right) ^{2
}-2\, \left(  \left( 2\,n \right) ! \right) ^{2}}{ \left(  \left( 2\,n
 \right) ! \right) ^{2}}},
$$
and its third through fifth standardized moments about the mean approach
$$
2\,{\frac {\sqrt {\pi } \left( \pi -3 \right) }{ \left( -\pi +4
 \right) ^{3/2}}},
 -{\frac {3\,{\pi }^{2}-32}{{\pi }^{2}-8\,\pi +16}},
4\,{\frac {\sqrt {\pi } \left( {\pi }^{2}+5\,\pi -25 \right) }{
 \left( -\pi +4 \right) ^{5/2}}}
$$
as $n \to \infty$.
\end{prop}
\subsection{Forward King walks}
Now we examine another special set of walks:
 \begin{defn}
 Let $\Knn$ be the set of walks from $(0,0)\to(n,n)$ with steps in $\{r,u,d\}=\{(1,0),(0,1),(1,1)\}$.  Let $K:=\bigcup_{n \in \N} \Knn$. For $w \in K$, let $n(w)$ be the $n$ such that $w \in \Knn$.
 \end{defn}
Think of $K$ as the set of journeys possible for a forward-marching King that end on the line $y=x$.  Each move, we take a step from $\{r,u,d\}$ (right, up, or diagonal).

For a set $E\subset K,$ define the generating function
$$
F_E(z,t):=\sum_{w\in E} z^{n(w)} t^{a_1(w)}.
$$

We we shall find an algebraic expression for $F_K(z,t)$.  The idea is to convert facts describing walks in $K$ to equations involving generating functions.

\begin{defn}
\leavevmode
\\
\begin{itemize}
\item Juxtaposition of two sets $A,B$ of walks denotes concatenation:
$$
AB:=\{w_1w_2 : w_1 \in A, w_2 \in B\}.
$$
If $A$ or $B$ is a singleton, we drop the braces: e.g., $aB:=\{a\}B$ for a walk $a$.
\item The \emph{Kleene star} of a set of walks is its closure under concatenation:
$$
E^*:=E \cup EE \cup \cdots = \{s_1s_2\cdots s_k : k \in \mathbb{N}, s_k \in E\}.
$$

\item We define the \emph{star of a generating function} $F$ to be
$$
F^*:=1+F+F^2+\cdots=\frac{1}{1-F}.
$$
\end{itemize}
\end{defn}

Now, let $N$ be the negative walks in $K$, i.e., walks satisfying $y<x$, save for the first and last points.  Let $\Phi$ be the nonpositive walks, i.e., walks in $y\leq x$.  Any negative walk is a right step followed by a nonpositive walk followed by an up step: $N=r\Phi u$.  Note that \emph{every} point of $w\in \Phi$ is counted as a losing time in $rwu$, so defining
$$
\tfp:=\sum_{w\in \Phi} z^{n(w)} t^{l(w)},
$$
we have
\begin{equation}
F_N=zt^2 \tfp.
\end{equation}\label{eq:fn}
Next, any nonpositive walk consists of diagonal steps and negative walks. So $\Phi=d^*(Nd^*)^*$,  and
\begin{equation}\label{eq:fp}
\tfp=(zt)^*(F_N(zt)^*)^*.
\end{equation} 

It is child's play for Maple to solve \eqref{eq:fn} and \eqref{eq:fp} for $F_N(z,t)$.  Now let $P$ be the set of positive walks, i.e., walks in $y>x$, except for the endpoints.  Positive walks are simply negative walks reflected about $y=x$, so $F_P=F_N(z,1)$ (all positive walks have zero losing times).

Finally, any forward King walk consists of diagonal steps, negative walks, and positive walks: $K=d^*(Nd^*\cup Pd^*)^*$.  In terms of generating functions,
$$
F_K=z^*(F_N z^*+F_P z^*)^*,
$$
and we are finished!  We have $F_K(z,t)$ as an algebraic expression.  Of course, it is rather messy, so we do not record it here.  To see it for yourself, use \verb+ForwardKingGF(z,t)+ in the Maple package. Unfortunately, $F_K$ is too complex to be amenable to either of the moment-finding methods discussed previously.  However, we should not lose hope\dots
\section{Numerically analyzing moment asymptotics}\label{sec:num}
\subsection{Recursive enumeration of  walks}
We started with steps in $\{(1,0),(0,1)\}$.  Then we added the diagonal step $(1,1)$.  Now let us take the affair even further.    

\begin{defn}
Given $S\subset \N^2$, let $W^S$ be the set of walks from $(0,0)$ with steps in $S$.  For $(a,b) \in \N^2$, let $W_{a,b}^S$ contain walks of $W^S$ ending at $(a,b)$.
\end{defn}

In the $S=\{(1,0),(0,1)\}$ case, we were able to calculate the moments of $a_1 | W_{n,n}^S$ in terms of $n$.  We cannot expect to do this in general.  Indeed, even in the (still very symmetric) case  $S=\{(1,0),(0,1),(1,1)\}$, we could not find nice expressions for the moments.  For even wilder choices of $S$, who knows if closed form generating functions over $W^S$ even exist?

However, we can \emph{fix} $(a,b) \in \N^2$ and focus on the \emph{finite} set of walks $W_{a,b}^S$.  Then the generating function 
$$
F_{a,b}(t):=\sum_{w \in W_{a,b}^S} t^{a_1}
$$
is a \emph{finite} polynomial in $t$, with easily computable moments.  Further, given fixed $S$ and $(a,b)$, we can make use of the fact that 
$$
W_{a,b}^S=\bigcup_{s \in S} W_{(a,b)-s}\{s\}
$$
to compute $F_{a,b}(t)$ with a recursive procedure; this is done in \verb+F2G+.  So for each $(a,b)$, the moments of $a_1 | W_{a,b}^S$ are easily computable.
\subsection{Asymptotic storybooks}
The procedure \verb+ChungFellerBook2D(S,M,K1,K2)+ uses \verb+F2G+ to compute the expectation, variance, and central standardized moments three through $M$ of $a_1 | W_{n,n}^{S'}$ for $S' \subset S,$ $n=K1,\dots,K2$.  It uses this data to guess the asymptotic behavior of the moments as functions of $n$.  We use the ansatzes $Cn$ for expectation, $Cn^2$, and $C$ for the third and higher central standardized moments.  

So, for each $S' \subset S$, a theorem about the asymptotic behavior of walks with steps in $S'$ is generated (step sets producing trivial theorems are automatically excluded).  

Of course, we must add the disclaimer that these ``theorems" are merely numerical approximations to the asymptotic behavior of the moments.  To be extra safe, the procedure \verb+ChungFellerBook2DSafe+ runs \verb+ChungFellerBook2D+ twice, with different $n$-ranges.  For each theorem it computes the constants twice; then it only keeps the agreeing digits.  

In the case $S=\{(1,0),(0,1)\}$, where we do know the moments as functions of $n$, we can confirm that \verb+ChungFellerBook2DSafe+ gives good results.

The table below summarizes the output of \\ \verb+ChungFellerBook2DSafe({[1,0],[0,1],[1,1],[2,0],[0,2]},6,100,110,+ \verb+190,200);+. 
The zeroth column is the set of allowed steps, where for brevity $ij:=[i,j]$.  Columns 1-6 are the expectation, variance, and limits of the third through sixth standardized central moments.  Note that by Proposition \ref{prop:cfm}, the exact values of the first row are  $n, n^2/3, 0, 1.8, 0, 27/7\approx .38571$.
\\
\\
\begin{tabular}{c c  c  c c c c} 
\hline
Steps & 1 & 2 & 3 & 4 & 5 & 6 
\\
\hline
$\{01, 10\}$ & $1.0000n$ & $0.3n^2$ & 0.0000 & 1.800 & 0.0000 & 3.86
\\
$\{01, 20\}$ & $0.38n$ & $0.1n^2$ & $0.0$ & $0.900$ & $-0.1$ & $1.93$
\\
$\{02, 20\}$ & $0.2500n$ & $0.043n^2$ & $0.0000$ & $0.900$ & $0.0000$ & $1.93$
\\
$\{01, 02,  10\}$ & $0.9n$ & $0.27n^2$ & $0.0$ & $1.8023$ & $-0.03$ & $3.87$
\\
$\{01, 02, 20\}$ & $0.33n$ & $0.07n^2$ & $-0.02$ & $0.90$ & $-0.1$ & $2.$
\\
$\{01, 10, 11\}$ & $0.8n$ & $0.2n^2$ & $0.0$ & $1.8$ & $0.$ & $3.9$
\\
$\{01, 11, 20\}$ & $0.666n$ & $0.15n^2$ & $0.001$ & $1.80$ & $0.01$ & $3.9$
\\
$\{02, 11, 20\}$ & $0.5n$ & $0.08n^2$ & $0.0$ & $1.80$ & $0.$ & $3.9$
\\
$\{01, 02, 10, 11\}$ & $0.81n$ & $0.22n^2$ & $0.$ & $1.80$ & $0.0$ & $3.9$
\\
$\{01, 02, 10, 20\}$ & $0.80n$ & $0.21n^2$ & $-0.01$ & $1.804$ & $0.$ & $4.$
\\
$\{01, 02, 11, 20\}$ & $0.6n$ & $0.1n^2$ & $-0.01$ & $1.80$ & $-0.1$ & $3.9$
\\
$\{01, 10, 11, 20\}$ & $0.81n$ & $0.22n^2$ & $0.$ & $1.80$ & $0.03$ & $3.89$
\\
$\{01, 02, 10, 11, 20\}$ & $0.75n$ & $0.19n^2$ & $-0.004$ & $1.8$ & $-0.011$ & $4.$
 \label{tab:output}
\end{tabular}
\subsection{Walks in higher dimensions}
This method easily generalizes to three or more dimensions.  If we consider walks in $\N^3$, then we have seven statistics to keep track of the number of times the walk visits the regions $x<y<z$, $x<z<y$, $y<x<z$, $y<z<x$, $z<x<y$, $z<y<x$, and ``none of the above."  The corresponding generating function (over walks to a fixed point in $\N^3$) is computed in \verb+F3G+.  Not surprisingly, this procedure is significantly slower than \verb+F2G+. 
\section{Conclusion}
Many areas still need to be explored.  For example, we have focused mainly on the number of losing times, $a_1$.  But the method of Section \ref{sec:num} could also be applied to the other statistics in Definition \ref{def:stat}.  Also, there is much to be done with walks in higher dimensions.  Our main goal in this paper was to illustrate by example that experimental mathematics can give us insights into problems where exact analysis is difficult.  Hopefully the results here are just the beginning. We encourage you to experiment with the Maple package and make discoveries of your own!
\section{Using the Maple package}
The Maple package accompanying this paper may be found at \url{http://www.math.rutgers.edu/~az202/Z/Feller.txt}.  To use it, place it in the working directory and execute \verb+read(`Feller.txt`);+.  To see the main procedures, execute \verb+Help();+. For help on a specific procedure, use \\ \verb+Help(procedure_name);+.  Happy exploring!
\section*{Acknowledgement} 
I thank Dr. Doron Zeilberger for introducing this project to me and  guiding my research in the right direction.
\bibliography{references}
\bibliographystyle{plain}

\end{document}